\documentclass{article}
\usepackage{amsmath}

\author{Gert Almkvist}
\title{POLYNOMIAL SOLUTIONS TO DIFFERENCE EQUATIONS CONNECTED TO PAINLEVE' II-VI}

\begin{document}

\maketitle

\section{\textbf{Introduction.}}

Around 1989 the following problem circulated: Given 
\[
\]
\[
a_0=a_1=a_2=a_3=1 
\]
and 
\[
a_na_{n-4}=a_{n-1}a_{n-3}+a_{n-2}^2 
\]
\[
\]
for $n\geq 4$. Show that all $a_n$ are integers (the $Somos$ $sequence$).
Soon various solutions appeared (see Gale [1] ).

The more general difference equation 
\[
\]
\[
a_na_{n-k}=\sum_{j=1}^{[k/2]}a_{n-j}a_{n-k+j} 
\]
\[
\]
where 
\[
a_0=a_1=\cdot \cdot \cdot =a_{k-1}=1 
\]
\[
\]
has integer solutions for $k=4,5,6,7$ but not for $k>7$. Why\ $7$ ? There
does not seem to be an explanation for this.

Here we shall instead consider similar nonlinear difference equations
involving polynomials connected to rational solutions of Painleve' equations
II-VI. We have the following table of equations 
\[
\]

P$_{II}:$ Yablonskii-Vorob'ev polynomials [2], [4]. 
\[
\]
\[
P_0=1\text{, }P_1=x. 
\]
\[
\]
\[
P_{n+1}P_{n-1}=-4(P_nP_n^{\prime \prime }-P_n^{\prime \text{ }2})+xP_n^2 
\]
\[
\]

P$_{III}:$ Umemura polynomials [3] 
\[
\]
\[
P_0=P_1=1
\]
\[
\]
\[
P_{n+1}P_{n-1}=-x^4(P_nP_n^{\prime \prime }-P_n^{\prime
}{}^2)-x^3P_nP_n^{\prime }+(cx+1)P_n^2
\]
\[
\]

P$_{IV}:$ [2]. 
\[
\]
\[
P_0=P_1=1 
\]
\[
\]
$^{}$%
\[
P_{n+1}P_{n-1}=P_nP_n^{\prime \prime }-P_n^{\prime \text{ }2\text{ }%
}+(x^2+2n-1)P_n^2 
\]
\[
\]
and 
\[
P_0=1,\text{ }P_1=x 
\]
\[
\]
\[
P_{n+1}P_{n-1}=P_nP_n^{\prime \prime }-P_n^{\prime \text{ }2}+(x^2+2n)P_n^2 
\]
\[
\]

P$_V:$ [5]. 
\[
\]
\[
P_0=P_1=1 
\]
\[
\]
\[
P_{n+1}P_{n-1}=x(P_nP_n^{\prime \prime }-P_n^{\prime \text{ }%
2})+P_nP_n^{\prime }+(\frac x8-v+\frac{3n}8)P_n^2 
\]
\[
\]

P$_{VI}:$ [6]. 
\[
P_0=P_1=1
\]
\[
\]
\[
P_{n+1}P_{n-1}=\frac{(x^2-4)^2}4(P_nP_n^{\prime \prime }-P_n^{\prime \text{ }%
2})+\frac{x(x^2-4)}4P_nP_n^{\prime }+(cx+(n-\frac 12)^2)P_n^2
\]
\[
\]
E.g. 
\[
y=\frac{P_n(1/x,c)P_{n+1}(1/x,c-1)}{P_n(1/x,c-1)P_{n+1}(1/x,c)}
\]
\[
\]
is a rational solution of P$_{III}$%
\[
\]
\[
y^{\prime \prime }=\frac{y^{\prime \text{ }2}\text{ }}y-\frac{y^{\prime }}x+%
\frac{ay^2+b}x+y^3-\frac 1y
\]
\[
\]
where (see [3] ) 
\[
a=2n-1+2c
\]
\[
b=2n+1-2c
\]
\[
\]
\[
\]

In this paper we consider the more general equation 
\[
\]
\[
P_{n+1}P_{n-1}=f(P_nP_n^{\prime \prime }-P_n^{\prime \text{ }%
2})+gP_nP_n^{\prime }+hP_n^2 
\]
with 
\[
P_0=P_1=1 
\]
\[
\]
We prove the following

\textbf{Theorem1 }Assume that $f,g,h$ are polynomials without common factor
and independent of $n$ . Assume also that consecutive $P_n$:s are without
common factor. Then all $P_n$ are polynomials if and only if 
\[
\]
\[
ff^{\prime \prime }-f^{\prime 2}+3f^{\prime }g-2fg^{\prime }-2g^2=0\text{
\qquad (*)}
\]

\[
\]

\textbf{Theorem2} The only solutions of (*) are the following 
\[
\]
(i) 
\[
f\text{ an arbitrary polynomial} 
\]
\[
g=\frac{f^{\prime }}2 
\]
\[
\]
(ii) 
\[
f=(ax+b)^k 
\]
\[
g=a(ax+b)^{k-1} 
\]
\[
\]
It turns out that these results also hold when 
\[
\]
\[
h=n+p(x) 
\]
\[
\]
which occurs in the P$_{IV}$ and P$_V$ cases. For the P$_{VI}$ case we have 
\[
\]
\[
h=n(n-1)+p(x) 
\]
\[
\]
then we get the condition 
\[
\]
\[
ff^{\prime \prime }-f^{\prime \text{ }2}+3f^{\prime }g-2fg^{\prime
}-2g^2+2f=0 
\]
\[
\]
which seems to have only solutions with $\deg (f)\leq 4$ . 
\[
\]
\[
\]
\textbf{2.The proofs.}

To shorten the proof we use the notation 
\[
P=P_n 
\]
\[
R=P_{n-1} 
\]
\[
S=P_{n-2} 
\]
\[
\]
We have to show that 
\[
\]
\[
f(PP^{\prime \prime }-P^{\prime \text{ }2})+gPP^{\prime }+hP^2 
\]
\[
\]
is divisible by $R$ . We have 
\[
\]
\[
SP=f(RR^{\prime \prime }-R^{\prime \text{ }2})+gRR^{\prime }+hR^2\equiv
-fR^{\prime \text{ }2}\text{ }mod\text{ }R 
\]
\[
\]
It follows 
\[
\]
\[
S^{\prime }P+SP^{\prime }=f^{\prime }(RR^{\prime \prime }-R^{\prime \text{ }%
2})+f(RR^{\prime \prime \prime }-R^{\prime }R^{\prime \prime })+g^{\prime
}RR^{\prime }+gR^{\prime \text{ }2}+gRR^{\prime \prime }+h^{\prime
}r^2+2hRR^{\prime } 
\]
\[
\]
\[
\]
and 
\[
S^2P^{\prime }\equiv R^{\prime \text{ }2}(S^{\prime }f-Sf^{\prime
}+Sg)-R^{\prime }R^{\prime \prime }Sf\text{ \qquad }mod\text{ }R 
\]
\[
\]
\[
\]
Differentiating once more we get 
\[
\]
\[
S^{\prime \prime }P+2S^{\prime }P^{\prime }+SP^{\prime \prime }\equiv
R^{\prime \text{ }2}(-f^{\prime \prime }+2g^{\prime }+2h)+R^{\prime
}R^{\prime \prime }(-2f^{\prime }+3g)-R^{\prime \prime \text{ }2}f 
\]
\[
\]
and 
\[
S^3P^{\prime \prime }\equiv R^{\prime \text{ }2}\left\{ S^2(-f^{\prime
\prime }+2g^{\prime }+2h)-S^{\prime \text{ }2}2f+SS^{\prime }(2f^{\prime
}-2g)+SS^{\prime \prime }f\right\} 
\]
\[
+R^{\prime \text{ }3}R^{\prime \prime }\left\{ S^2(-2f^{\prime
}+3g)+SS^{\prime }2f\right\} +R^{\prime \prime \text{ }2}\left\{
-S^2f\right\} 
\]
\[
\]
Hence 
\[
S^4(PP^{\prime \prime }-P^{\prime \text{ }2})\equiv R^{\prime \text{ }%
4}\left\{ S^2(ff^{\prime \prime }-2fg^{\prime }+2f^{\prime }g-f^{\prime 
\text{ }2}-2fh-g^2)+S^{\prime \text{ }2}f^2-SS^{\prime \prime }f^2\right\} 
\]
\[
-R^{\prime \text{ }3}R^{\prime \prime }S^2fg 
\]
\[
\]
and 
\[
S^4PP^{\prime }\equiv R^{\prime \text{ }4}\left\{ S^2(ff^{\prime
}-fg)-SS^{\prime }f^2\right\} +R^{\prime \text{ }3}R^{\prime \prime }S^2f^2 
\]
\[
\]
\[
S^4P^2\equiv R^{\prime \text{ }4}S^2f^2 
\]
\[
\]
Adding up we obtain 
\[
\]

\[
S^4\left\{ f(PP^{\prime \prime }-P^{\prime \text{ }2})+gPP^{\prime
}+hP^2\right\} \equiv 
\]
\[
R^{\prime \text{ }4}f\left\{ S^2(ff^{\prime \prime }-f^{\prime \text{ }%
2}+3f^{\prime }g-2fg^{\prime }-2g^2)-f(f(SS^{\prime \prime }-S^{\prime \text{
}2})+gSS^{\prime }+hS^2)\right\} \equiv 0\text{ \qquad }mod\text{ }R 
\]
\[
\]
if and only if 
\[
ff^{\prime \prime }-f^{\prime \text{ }2}+3f^{\prime }g-2fg^{\prime }-2g^2=0 
\]
\[
\]
since 
\[
f(SS^{\prime \prime }-S^{\prime \text{ }2})+gSS^{\prime }+hS^2 
\]
\[
\]
is divisible by $R.$%
\[
\]

We also use that $(R,S)=1$ and $(R,R^{\prime })=1$ . If $(R,R^{\prime })\neq
1$ then $R$ would have a double root causing $(P,R)\neq 1$ considered in $%
C[x].$%
\[
\]

\textbf{Remark: }In the cases P$_{II}$ and P$_{III}$ one can assume that $P_n
$ and $P_{n+1}$ have no common roots so the assumptions in Theorem 1 are
fulfilled. But experiments show that Theorem 1 is true without any assum

ptions on the $P_n$:s.
\[
\]

\textbf{Proof of Theorem 2:}

One checks that 
\[
g=\frac{f^{\prime }}2 
\]
\[
\]
solves 
\[
ff^{\prime \prime }-f^{\prime \text{ }2}+3f^{\prime }g-2fg^{\prime }-2g^2=0 
\]
\[
\]
To find all solutions we make the substitution 
\[
g=u+\frac{f^{\prime }}2 
\]
\[
\]
Then we get 
\[
fu^{\prime }-\frac{f^{\prime }}2u+u^2=0 
\]
\[
\]
Put 
\[
w=\frac 1u 
\]
\[
\]
Then 
\[
\sqrt{f}w^{\prime }+\frac w{2\sqrt{f}}=(\sqrt{f}w)^{\prime }=\frac 1{\sqrt{f}%
} 
\]
\[
\]
with the solution 
\[
w=\frac 1{\sqrt{f}}\left\{ \int \frac{dx}{\sqrt{f}}+C\right\} 
\]
\[
\]
and 
\[
u=\frac{\sqrt{f}}{\int \frac{dx}{\sqrt{f}}+C} 
\]
\[
\]
The only way $u$ can become a polynomial is when $C=0$ and 
\[
\]
\[
f=(ax+b)^k 
\]
\[
\]
where $k>0$ is an integer. Then 
\[
\int \frac{dx}{\sqrt{f}}=\int (ax+b)^{-k/2}dx=\frac{(ax+b)^{1-k/2}}{a(1-k/2)}
\]
\[
\]
and 
\[
u=a(1-\frac k2)(ax+b)^{k-1} 
\]
\[
\]
Finally 
\[
g=u+\frac{f^{\prime }}2=a(ax+b)^{k-1} 
\]
\[
\]
\[
\]

If 
\[
h=h_n=n+p(x) 
\]
\[
\]
we will end up with the following coefficient of $S^2$%
\[
\]
\[
-2h_{n-1}+h_n=-2(n-1+p(x))+n+p(x)=-(n-2+p(x)=-h_{n-2} 
\]
\[
\]
so nothing is changed in the proof of Theorem 1. 
\[
\]
\[
\]

If 
\[
h=h_n=n(n-1)+p(x) 
\]
\[
\]
then 
\[
-2h_{n-1}+h_n=-2((n-1)(n-2)+p(x))+n(n-1)+p(x)= 
\]
\[
-((n-2)(n-3)+p(x))+2=-h_{n-2}+2 
\]
\[
\]
which means that the condition (*) is changed to 
\[
\]
\[
ff^{\prime \prime }-f^{\prime \text{ }2}+3f^{\prime }g-2fg^{\prime
}-2g^2+2=0 
\]
\[
\]
If you try to find polynomial solutions with $\deg (f)>4$ you will run into
contradictions.E.g. there are the following solutions 
\[
\]

(i)

$f=x^2+ax+b$

$g=2x+a$

\[
\]

(ii)

$f=(x^2-1)^2$

$g=x(x^2-1)$%
\[
\]

(iii)

$f=(x^2+2x)^2$

$g=(x+1)(x^2+2x)$%
\[
\]

(iv)

$f=\dfrac{(x^2-4)^2}4$

$g=\dfrac{x(x^2-4)}4$

\[
\]

\textbf{References:}

1.D.Gale, The strange and surprising saga of the Somos sequences, Math.
Intelligencer 13 (1991), no 3, 40-42 and no 4, 49-50. 

2.S.Fukutani, K.Okamoto, H.Umemura, Special polynomials and the Hirota
bilinear relations of the second and the fourth Painleve' equations, Nagoya
Math. J. 159 (2000), 179-200. 

3.K.Kajiwara, T. Masuda, On the Umemura polynomials for the Painleve' III
equation, Physics Letters A 260 (1999), 462-467 

4.M.Kaneko, H.Ochiau, On the coefficients of Yablonskii-Vorob'ev
polynomials, QA 0205178 

5.M.Noumi, Y.Yamada, Umemura polynomials for the Painleve' V equation,
Physics Letters A 247 (1998), 65-69 

6.M.Taneda, Polynomials associated with an algebraic solution of the sixth
Painleve' equation, Japan. J. Math. (N.S.) 27 (2001), 257-274.

\end{document}